\documentclass[12pt]{article}
\usepackage{mathpazo}
\usepackage{amsfonts,amssymb,amsmath}
\usepackage{enumerate}
\usepackage{epigraph}
\usepackage{color,soul}     
\usepackage[normalem]{ulem}  

\newtheorem{proposition}{Proposition}

\newtheorem{example}{Example}
\newtheorem{definition}{Definition}

\newenvironment{proof}[1][Proof]{\noindent\textbf{#1.} }{\ \rule{0.5em}{0.5em}}

\begin{document}

\title{\textbf{Boolean-Valued Sets as Arbitrary Objects}}
\author{Leon Horsten \\ University of Konstanz}
\maketitle

\begin{abstract}
\noindent This article explores the connection between boolean-valued class mo\-dels of set theory and the theory of arbitrary objects in roughly Kit Fine's sense of the word. In particular, it explores the hypothesis that the set theoretic universe as a whole can be seen as an arbitrary entity, which can in turn be taken to consist of arbitrary objects, viz. arbitrary sets.
\end{abstract}



\section{Introduction}

Contemporary philosophy of physics aims to develop metaphysical interpretations of fundamental current physical theories. In philosophy of quantum mechanics, for instance, researchers articulate metaphysical accounts of what the physical world at the micro-level \emph{could be} like given our current quantum mechanical theories.
The aim of this article is to do something similar for set theory. The aim is  to articulate a new metaphysical view of what the set theoretic world could be like given our current set theoretic theories and practices.

The most important development in set theory since the second world war is Cohen's discovery of \emph{forcing}, which is an incredibly powerful and flexible technique for producing independence results.
There are today two main approaches to forcing. The first approach is sometimes called the \emph{forcing poset approach}.\footnote{See \cite{Kunen1980}.} The second approach is called the \emph{boolean-valued approach}. The latter was pioneered by Scott and Solovay (and independently discovered by Vop\u{e}nka), and was first described in \cite{Scott1967}.\footnote{The early history of forcing is described in \cite{Moore1988}.}
The two approaches are for most intents and purposes equivalent.

The boolean-valued approach, as described in \cite{Bell2005}, is centred around the concept of boolean-valued sets, which are functions into a complete boolean algebra. Boolean-valued sets have been studied mostly with the aim of proving set theoretic independence results. Here I want to consider structures of boolean-valued sets from a metaphysical perspective. I will argue that boolean-valued sets can be seen as \emph{arbitrary objects} in the sense of \cite{Fine1985} and \cite{Horsten2019}. Indeed, Fine himself suggested that arbitrary object theory might be applicable to forcing \cite[p.~45--46]{Fine1985}, although his suggestion has hitherto not been followed up.

I will develop the metaphysical hypothesis that there is a sense in which the set theoretic universe itself is also an arbitrary entity. On the view that I explore, there is only one mathematical universe. But just as the elements in it, the set theoretic universe as a whole is an arbitrary entity. And just as the arbitrary sets in the universe can be in different states, the set theoretic universe can also be in many mutually incompatible states.

In this article the boolean algebra-approach to forcing is used as a tool to express a metaphysical view.
Given the mathematical equivalence of the boolean algebra-approach and the partial order-approach (``poset-approach'') to set forcing, the metaphysical view that I want to explore can also be expressed using the partial order-approach, but I will not do so here. For \emph{class forcing}, the two approaches are not mathematically equivalent.\footnote{See \cite{Antosetalforthc}.} In this article, I will mostly ignore class forcing.

The proposal that is explored in this article is tentative, and is deliberately kept ``open'' at several junctures. That is, there are multiple ways in which it can be developed more fully. This is just because at this point I am unsure what the best way is of fleshing out my proposal in a detailed way. Also, I will not attempt to argue that it is, all things considered, more plausible than rival proposals: that task is left for a future occasion.

The structure of this article is as follows. First, I review key elements of the theory of arbitrary objects. Then I will show how in an obvious way arbitrary objects play a role in the model that is described in \cite{Scott1967}, which is the oldest incarnation of the boolean-valued models approach to forcing. Then I will show how arbitrary object theory can also be used to give a metaphysical interpretation of `modern' boolean-valued models.

In what follows I use the notation that is used in \cite{Bell2005}. In particular, concerning algebraic notions, I denote the join, meet, complementation, and implication operations as $\vee$, $\wedge$, $^c$, and $\Rightarrow$, respectively, and I denote the top and bottom elements of an algebra as 1 and 0, respectively. Moreover, I use the slightly unusual convention by Bell to refer to boolean algebras as $B$, $B'$,\ldots


\section{Arbitrary objects}\label{arbitrary objects}

An \emph{arbitrary} $F$ is an abstract object that can be in a \emph{state} of being some or other $F$. We may say that an arbitrary $F$ \emph{coincides} with some $F$ in a state, or takes a certain \emph{value} in some state.  So, mathematically, 
an arbitrary $F$ can be modelled as a function $$ f: \Omega \rightarrow F, $$ where $\Omega$ is a \emph{state space}, and $F$ is a collection of objects. In order to develop a basic feeling for what arbitrary objects are like, let us briefly consider a few simple examples.

\begin{example}
Consider \emph{the man on the Clapham omnibus}. Such an arbitrary object could (in some sense) be me, or it could be my next door neighbour. But this arbitrary object is neither identical with me, nor with my next door neighbour.
\end{example}

\begin{example}\label{example arbitr nat numb}
Consider an \emph{arbitrary natural number}. Such an arbitrary object can be in the state of being the number 3, but it can also be in the state of being the number 4.
\end{example}

\noindent Let us briefly relate this second example back to the modelling proposal above.\footnote{For more details, see \cite[section 4.1.2]{Horsten2019}.} In example \ref{example arbitr nat numb}, we evidently take the value space to be $\mathbb{N}$. We want there to be as many states as are needed in order to have fully arbitrary natural numbers, i.e. arbitrary natural numbers that can be any specific natural number. Moreover, we do not seem to have reasons to have more states.  This means that $\lvert  \Omega \rvert= \aleph_0$. Arbitrary natural numbers can then be seen as threads through the matrix $\mathbb{N} \times \mathbb{N}$, i.e. functions from $\mathbb{N}$ to $\mathbb{N}$.\label{arbitrary nat numb}

Typically, for a property $F$, there are more than one arbitrary $F$'s.\footnote{Fine holds that for every $F$, there is ultimately no more than one ``independent'' arbitrary $F$ \cite[p.~69]{Fine1983}. I will not make this assumption here.} For instance, consider an arbitrary natural number $a_1$ strictly between 3 and 6. Then there is also another arbitrary natural number strictly between 3 and 6, call it $a_2$, which in every state differs from $a_1$. So, for instance, in a state where $a_1$ takes the value 4, $a_2$ takes the value 5. This shows that arbitrary $F$'s can be \emph{correlated} with each other.\footnote{Note that this means that it is strictly speaking wrong to speak of \emph{the} man on the Clapham omnibus.}

It has been argued, a.o. by Frege,\footnote{See \cite{Frege1904}.} that there are no arbitrary objects, and this is still the prevailing view. But in the spirit of \cite{Fine1985} and \cite{Horsten2019}, I will take arbitrary objects ontologically seriously. The aim of this article is not to argue for this metaphysical stance.

In many cases, the function range of an arbitrary $F$, when regarded as a function, consists of \emph{specific} objects. For instance, in a state where $a_1$ coincides with the number 4, it takes a \emph{specific} value. But there are also arbitrary objects that can be in a state of being this or that \emph{arbitrary} object. For instance, an \emph{arbitrary} arbitrary natural number strictly between 3 and 6 can be in a state of being the arbitrary number $a_1$, but it can also being in a state of being the arbitrary number $a_2$. Such higher-order arbitrariness will play an important role in what follows.


I will also be liberal in not just considering maximally specific state descriptions (also known as \emph{possible worlds}). Instead, I will also permit as states situations that are less than fully specified: call them \emph{partial states}. These partial states can then be modelled as \emph{sets} of possible worlds.


\section{Forcing and random variables}

In early work on boolean-valued models, \emph{random variables} play an important role. In particular, this is so in the first exposition of the method of boolean-valued models, Scott's beautiful article \emph{A proof of the independence of the Continuum Hypothesis} \cite{Scott1967}.

Scott starts his construction of boolean-valued models with a probability triple  $\langle \Omega, \mathcal{A}, \mathrm{P} \rangle$, where $\Omega$ is a state space, $\mathcal{A}$ is a $\sigma$-algebra on $\Omega$, and $\mathrm{P}$ is a probability function defined on $\mathcal{A}$. This probability triple is the background of the notion of a \emph{random real} over $\Omega$, where a random real over $\Omega$ is a function
$$ \xi : \Omega \rightarrow \mathbb{R}   $$
that satisfies some measurability constraint.\footnote{In particular, it is required that for each $r\in \mathbb{R}:$
$$ \{  o \in \Omega: \xi(o) \leq r  \}  $$
is measurable.}
Let $\mathcal{R}$ be the collection of random reals. It is easy to see that $\mathbb{R}$ is canonically embedded in $\mathcal{R}$ (by constant functions).

Scott's aim is, roughly, to construct a boolean-valued analogue of the classical rank $V_{\omega + 2}$, which is the level of the iterative hierarchy where the Continuum Hypothesis ($\mathrm{CH}$) is decided. In this boolean-valued model the axioms of set theory insofar as they describe $V_{\omega + 2}$, turn out to be true, whereas $\mathrm{CH}$ is false in this model.

The language in which Scott describes the initial transfinite levels of the iterative hierarchy has a type-theoretic flavour.\footnote{But this is not essential for his argument, as Scott himself observes.} In particular, it contains variables ranging over real numbers, and variables ranging over functions on the reals. The set of natural numbers $\mathbb{N}$ is \emph{defined} in this language as a special collection of reals \cite[p.~95]{Scott1967}.

In the resulting model $\mathcal{S}$, the real number variables range over random reals (as defined above). The function variables range over a set $\mathcal{R}^{\mathcal{R}}$ of functions from $\mathcal{R}$ to $\mathcal{R}$ that meet an extensionality condition.\footnote{See \cite[p.~102]{Scott1967}.}

In the boolean-valued model $\mathcal{S}$, sentences of the language take values in a complete boolean algebra $B$, which is obtained from the boolean $\sigma$-algebra $\mathcal{A}$ by identifying events that differ from each other only by a set of probability 0 (as measured by the probability function $\mathrm{P}$).
Moreover, $ B $ can be seen to have the countable chain condition, which entails that $ B $ is \emph{complete}.

Then Scott chooses $\Omega$ in such a way that $\mathcal{S}$ contains \emph{many} random reals that are ``orthogonal'' to each other. This ensures that $\mathcal{S} \models \neg \mathrm{CH}$, where $\models$ is the \emph{boolean-valued} truth relation. In particular, the ``degree'' to which two random reals $\xi$ and $\eta$ coincide according to $\mathcal{S}$ is ``measured''  by a boolean value, i.e., an element of $B$. And such an element of $B$ can \emph{roughly}\footnote{I.e., up to $P=0$ difference.} be taken to be the set of states on which $\xi$ and $\eta$ coincide. 

Now $\mathcal{S}$ only verifies the usual set theoretic axioms as far as $V_{\omega + 2}$ goes.  But Scott sketches how $\mathcal{S}$ can fairly routinely be extended to a boolean-valued model of $ZFC$ that still makes $\mathrm{CH}$ false.

As objects that take values in states, Scott's random reals are \emph{arbitrary objects} in the sense of \cite{Fine1985} and \cite{Horsten2019}\footnote{See in particular \cite[chapter 10]{Horsten2019}.} (or at least they are modelled in the same way). But the values of \emph{function variables} are \emph{not} natural modellings of arbitrary objects. Going up the hierarchy, functionals, etcetera, are also \emph{not} arbitrary objects. This ``non-uniformity'' is eliminated in later versions of boolean-valued model theory, such as \cite{Bell2005}, as we will see shortly.

The take-away message is that arbitrary objects have played a role in boolean-valued models from the start. Random variables in Scott's sense have mostly disappeared from modern treatments of boolean-valued models,\footnote{But not entirely: see for instance \cite{Krajicek2011}.} and Scott himself already observed that his method for proving the independence of the Continuum Hypothesis does not really require them \cite[p.~110]{Scott1967}. 
Nevertheless, I will argue that in more recent versions of boolean-valued model theory, arbitrary objects play an even more pervasive role, albeit in a somewhat less obvious way.



\section{Boolean-valued sets as arbitrary objects}\label{boolean valued sets}

Let us now turn to the contemporary approach to boolean-valued models, as described in \cite{Bell2005}.

The domain of a boolean-valued class model $V^{(B)}$ consists of functions 
$$ u:  V^{(B)} \rightarrow B ,$$ where $B$ is a complete Boolean algebra. $dom(u)$ can be seen as the \emph{quasi-elements} of $u$. And the elements of $dom(u)$ are \emph{themselves} boolean-valued sets \cite[p.~21]{Bell2005}. This is reflected in the recursive build-up of the universe $V^{(B)}$ of boolean-valued sets.

Given the Stone representation theorem, the boolean algebra $B$ can be conceived of as a field of sets.\footnote{Uniquely so if $B$ is atomic. Note that the Stone space associated with a complete boolean algebra is not necessarily a complete algebra.} Each element of $B$ can then be seen as a \emph{set of possible worlds}, i.e., a partial (or total) state. In other words, $B$ is an algebra of states, where the join operation expresses union of states (`$a$ or $b$'). If in $B$ we have $a<b$, then the state $a$ is a \emph{refinement} of state $b$. The algebra $B$ need not be atomic: there may be no maximally specific states (`state descriptions' in the Carnapian sense). \emph{Partitions of unity} (as particular anti-chains in $B$) are then especially significant as collections of mutually exclusive and jointly exhaustive sets of states. In the absence of atomicity, partitions of unity are the  the closest counterparts to the set of all Carnapian possible worlds, which are possible worlds that are maximally specific and therefore decide all propositions.

Actuality plays no role in the picture. Just as it makes no sense to ask which state \emph{the fair coin} (an arbitrary object!) is actually in (heads or tails), there is no state that  $V^{(B)}$ is actually in.\footnote{Similarly, a potentialist about the set theoretic universe should hold that it makes no sense to ask which sets \emph{actually} exist.} There are just many states that  $V^{(B)}$ \emph{can} be in. Certainly the maximally unspecific top element $1 \in B$ should not be seen as the actual world. If there are no atoms in $B$, then there is not even a \emph{candidate} of being the actual world in the Carnapian sense.

So far, we have only defined the domain of a boolean-valued class model. The definition of a boolean-valued class model $V^{(B)}$ is incomplete without a specification of the interpretation of logical and non-logical vocabulary. Like Bell, we abbreviate $[\![ \varphi ]\!]^{V^B}$ as $[\![ \varphi ]\!]^{B}$, and we  consider identity and elementhood in some more detail. We define \cite[p.~23, 1.15]{Bell2005}:
\begin{equation}\label{=}  [\![ u \in v ]\!]^B  \equiv \bigvee_{y \in dom (v)} (v(y) \wedge [\![ y=u ]\!]^B) .\end{equation}
\noindent This is what it means for some boolean-valued set $u$ to be \emph{to some extent} a member of the boolean-valued set $v$. The extent is measured by a boolean value.

Given extensionality, identity and elementhood are intertwined in set theory: identity also constitutively depends on elementhood. So in the boolean-valued framework we have \cite[p.~23, 1.16]{Bell2005}:

\begin{equation}\label{in}
\begin{split} [\![ u = v ]\!]^B  \equiv  \hspace{4cm} \\  \bigwedge_{y \in dom (v)} ( v(y) \Rightarrow  [\![ y \in u]\!]^B) \wedge \bigwedge_{y \in dom (u)} (u(y) \Rightarrow  [\![ y \in v]\!]^B)    ,
\end{split}
\end{equation}

\noindent where $x \Rightarrow y$ is an abbreviation of $x^c \vee y$ (with $^c$ being the complementation operation of $B$). This completes the definition of the notion of a boolean-valued class model $V^{(B)}$. Note that these models are in general not classical, in the sense that they take their values in a boolean algebra that need not be the two-element boolean algebra---although they do make all laws of classical logic true.

The boolean-valued truth conditions of non-atomic statements are exactly what you would expect \cite[p.~22]{Bell2005}, so there is no need to spell them out here. 
These compositional clauses determine a notion of boolean-valued truth. From now on I will write $V^{(B)} \models \phi$ for $[\![ \phi ]\!]^B = 1$.

Now I want to argue that \emph{all} elements of $V^{(B)}$ can straightforwardly be seen as arbitrary objects. Thus the non-uniformity of Scott's model, where the range of the first-order quantifiers is somehow distinguished, is eliminated.

We have seen that $B$ can be seen as a state space. But given that boolean-valued sets are functions $ u:  V^{(B)} \rightarrow B $, boolean-valued sets are arrows that ``point in the wrong direction'' for being arbitrary objects. Their domain, rather than their range, should be a state space. 

But we will see that this problem can easily be remedied.
First I explain how $V^{(B)}$ \emph{itself} can be seen as an arbitrary entity. Then I describe how the elements of $V^{(B)}$ can also be seen as arbitrary objects.

A `partial' state can be seen as a \emph{situation}.\footnote{This is inspired by Barwise and Perry, who also call partial possible worlds \emph{situations}. However, they model situations in a very different way: see, e.g., \cite{BarwisePerry1983}.} Then a maximal anti-chain is a collection of mutually exclusive but jointly exhaustive situations.
 Every situation $a$ can, as we will see shortly, itself be seen as a boolean-valued universe $V^{(B)}_a$. In particular, $V^{(B)}_a$ will make all the principles of $ZFC$ true. So if $\{ a_1,\hdots , a_k, \hdots  \}$ is a maximal anti-chain (partition of unity) in $V^{(B)}$,  then $V^{(B)}$  can be in the state of being $V^{(B)}_{a_k}$. More precisely: in the situation $a_k$, the set theoretic universe is in the state of being $V^{(B)}_{a_k}$.

Consider a boolean-valued set $u \in V^{(B)}$. Then for any situation $a_k$, $u$ ``is'' a boolean-valued set $u_{a_k}$ at $a_k$, where $u_{a_k}$ is defined as follows:

\begin{center}
\item For all $s \in dom(u)$: $u_{a_k}(s) \equiv a_k \wedge u(s) .$
\end{center}
Thus boolean-valued sets can be seen as \emph{higher-order} arbitrary objects in the sense of \cite[section 6.9]{Horsten2019}: they are objects that can be in states that are \emph{themselves} arbitrary objects.



Now we can say exactly what $V^{(B)}_a$ is, for any $a \in B$: $V^{(B)}_a$ is the boolean-valued class model generated by the restricted boolean algebra $B_a$.\footnote{The restricted boolean algebra $B_a$ consists of all elements of the form $y\wedge a$, with $y\in B$, and which is such that for $x,y \in B_a$, $x \star y$ is the same as  $x \star y$ in $B$ for $\star \in \{ \wedge, \vee  \}$, and $x^c$ in $B_a$ is $x^c \wedge a$ in $B$ \cite[p.~79]{Jech2006}.}  It is then easy to see that $V^{(B)}_a = V^{(B_a)}$.

If $B$ is a complete boolean algebra, so is the restricted algebra $B_a$. If $B$ has the countable chain condition, then $B_a$ has it also. And observe that for all $a \in B$, $V^{(B)}_a \models ZFC$. On the other hand, the structure $V^{(B)}$ might be such that neither $\mathrm{CH}$ nor $\neg \mathrm{CH}$ is true in it, but that it \emph{could be} in a state where $\mathrm{CH}$ is true and it \emph{could be} in a state where $\mathrm{CH}$ is false. Such a $V^{(B)}$ could function as a toy model of a set theoretic universe in which neither $\mathrm{CH}$ nor $\neg \mathrm{CH}$ is true.

Let us now turn to the problem of ``reversing the arrows''.

\begin{definition}\label{u star}
For any $u \in V^{(B)}$, $u^*$ is the function such $u^*(a) = u_a$ for every $a \in B$.
\end{definition} 

\begin{proposition}

\hspace{0cm}

\begin{enumerate}
\item $u^*$ is uniquely determined by $u$;
\item $u$ is uniquely determined by $u^*$.
\end{enumerate}

\begin{proof}
\noindent Clause 1 follows from definition \ref{u star}.

\noindent Clause 2. of the proposition follows because $u = u_1$ (where $1$ is, as before, the top element of $B$).
\end{proof}

\end{proposition}

So whether we take $V^{(B)}$ to consists of boolean-valued sets $u$ or their counterparts $u^*$ makes no mathematical difference.
But the $u^*$'s are functions from the state space $B$ to $V^{(B)}$. Therefore they are \emph{arbitrary objects}, or, to be philosophically more correct, they are natural representations of arbitrary objects. Thus we can regard every $u\in V^{(B)}$ as an arbitrary object in the sense of section \ref{arbitrary objects}. More in particular:
\begin{itemize}
\item the $u^*$'s are \emph{total} arbitrary objects (since $dom(u^*) = B$);
\item the state space $B$ of the $u^*$'s consists mostly of \emph{partial} states (since typically the algebra $B$ will be non-atomic).
\end{itemize}

In category theoretic terms, what we are doing is treating Boolean-valued sets as \emph{variable sets}. Given a boolean algebra $B$, one can consider the category $\mathrm{Shv}(B)$ of sheaves over $B$. Then it can be shown that $\mathrm{Shv}(B)$ and $V^{(B)}$ are \emph{equivalent} as categories \cite[p.~180]{Bell2005}.\footnote{The connection between sheaf theory and arbitrary object theory was already observed by Fine: see \cite[p.~47]{Fine1985}.}

Thus both $V^{(B)}$ and boolean-valued sets $u$ in $V^{(B)}$ are arbitrary objects in the following sense. If the boolean algebra $B$ is atomless, then as we ``go down'' $B$, the universe $V^{(B)}$ takes a more specific state $V^{(B)}_a$, without ever reaching a maximally specific state. Likewise, as we go down $B$, a typical boolean-valued set $u$ takes a more specific state $u_a$, without ever reaching a maximally specific state.



\section{Kaleidoscopic absolutism}

We have seen how boolean-valued universes, and the sets that they contain, can be seen as arbitrary objects. Now I will argue that the set theoretic universe as a whole can itself be seen as an arbitrary entity. The slogan is, roughly:
\begin{quote}
The set theoretic universe is the \emph{arbitrary} $V^{(B)}$. 
\end{quote}
Let us designate the arbitrary $V^{(B)}$ as $\mathcal{V}$. 

As in the case of arbitrary natural numbers (p.~\pageref{arbitrary nat numb}), $\mathcal{V}$ can be modelled as a function from a state space to a value space. The value-range of $\mathcal{V}$ is a collection of $V^{(B)}$'s where $B$ is an element of a \emph{collection} $\mathcal{B}$ of complete boolean algebras. Since we do not need more states than there are elements of $\mathcal{B}$, the collection $\mathcal{B}$ can then be taken to be the state space of $\mathcal{V}$. So $\mathcal{V}$ can be modelled as the (possibly proper class-size) function that takes each $B \in \mathcal{B}$ to $V^{(B)}$.

At this point, my account becomes somewhat vague: I am not able to say with much precision what $\mathcal{B}$ and $V$ are like. $\mathcal{B}$ should be a \emph{large} collection of complete boolean algebras (possibly proper class-sized), so that what set theorists regard as real possibilities are all represented as states that $\mathcal{V}$ can be in. $V$ should be large enough  to maximise the interpretive power of set theory \cite[section 5]{Steel2014}. Beyond this, I see only a few more constraints that $\mathcal{B}$ and $V$ satisfy: see p.~\pageref{constraints} below. 

You may ask: can we not ``complete'' the state space of $\mathcal{V}$ to a \emph{complete boolean algebra} $\mathcal{B}^*$ and take the set theoretic universe to be (structurally like) $V^{(\mathcal{B}^*)}$? But this does not work. As is pointed out for instance in \cite{Antosetalforthc}, if $\mathcal{B}$ is indeed a proper class, then  $\mathcal{B}^*$ is a hyperclass, and $V^{(\mathcal{B}^*)}$ therefore does not make $ZFC$ true. 

So on the proposed view, the set theoretic universe is an arbitrary entity that \emph{can} be structurally like a $V^{(B)}$, and only like some $V^{(B)}$; but this arbitrary entity $\mathcal{V}$ itself \emph{is not} structurally like a $V^{(B)}$. Here the requirement that the set theoretic universe can \emph{only} be structurally like a $V^{(B)}$ is an expression of confidence of many set theorists that the only way of generating set theoretic independence---aside from the G\"odelian independence phenomena---is by forcing techniques. It is of course not guaranteed that this confidence will in the long run turn out to be justified.

Like all slogans, the motto that the set theoretic universe is the \emph{arbitrary} $V^{(B)}$ has to be taken with a grain of salt. Because of well-known anti-reductionist arguments,\footnote{See for instance \cite{Benacerraf1965}.} the thesis should not be that the set theoretic universe is an entity that can be in the state of \emph{being} this or that $V^{(B)}$. After all, just as it is unreasonable to hold that the number 19 \emph{is} some pure set or other, so it is unreasonable to maintain that the set theoretic universe can be in a state of being some $V^{(B)}$. The point is rather that it can be in  states that are structurally like, or can be fruitfully \emph{modelled} as, $V^{(B)}$'s.

It is sometimes argued that there are different, equally valid concepts of set, and that it is somehow indeterminate which of these notions is described in set theory.\footnote{See for instance \cite[p.~416]{Hamkins2012}, but also \cite{NodelmanZalta2014}.} The position that I am putting forward here is not intended as an articulation of this view. The thought that I am trying to develop is \emph{not} that two states of the set theoretic universe describe different set \emph{concepts}. Rather, the view is that there is \emph{one} conception of set  that the set theoretic universe answers to: a notion of set as an arbitrary object. 

The central component of the proposed view consists of truth definitions for the formulas of the language of set theory ($\mathcal{L}_{ZFC}$).
I have sketched the definition of truth in a boolean-valued structure in section \ref{boolean valued sets}. But a natural definition of truth in the set theoretic universe $\mathcal{V}$ can also be given:
$$ \mathcal{V} \models \varphi = : \textrm{ for all } V^{(B)} \textrm{ in the value range of } \mathcal{V}:  V^{(B)} \models \varphi . $$
\noindent As adumbrated above, this truth definition is somewhat vague: we do not have a strong grasp of what the range of $ \mathcal{V}$ is.

By Tarski's theorem on the undefinability of truth, this definition for truth in $ \mathcal{V} $ can only be expressed in an extension $\mathcal{L}_{ZFC}^+$ of the language $\mathcal{L}_{ZFC}$. This `boolean-valued' truth definition quantifies over boolean-valued models ($V^{(B)}$'s), which consist of boolean-valued sets  that are applied to other boolean-valued sets, boo\-lean operations are applied to values of boolean-valued sets, and so on. We have seen that every boolean-valued set can be seen as an \emph{arbitrary object}, and that also the $V^{(B)}$'s themselves can in the same sense be seen as arbitrary objects. So we could re-write the boolean-valued truth definition in terms of arbitrary objects instead of boolean-valued sets---although I will not do so here. Thus, ultimately, the view that I am considering presupposes an independent grasp of arbitrary objects, states, and values of arbitrary objects. In this sense, the boolean-valued truth definition is intended to be interpreted \emph{in arbitrary object theory}, which is a part of metaphysics. The idea is that $\mathcal{L}_{ZFC}^+$ is \emph{itself} interpreted as being about arbitrary entities.

The view that is proposed is meant to be a \emph{foundational} interpretation. As such, it must stand on its own two legs. It must not be parasitic on any other foundational interpretations, in particular standard interpretations of set theory exclusively in terms of ``$\{ 0,1 \}$-valued sets''. The proposed view should be \emph{logically and conceptually autonomous} from them, in the sense that it should be possible to \emph{state} it without appealing to notions belonging specifically to such interpretations, and that it should be possible to \emph{understand} it without first understanding the notion of $\{ 0,1 \}$-valued set. \cite[p.~241]{LinneboPettigrew2011}.

You might worry that this autonomy requirement is not satisfied: the boolean-valued sets (and hence also the arbitrary objects) to which the account appeals are defined (recursively) in terms of quantification over $V$.
But this concern is misguided. The boolean-valued truth definition quantifies over the $V^{(B)}$'s and their elements, both of which are interpreted as arbitrary objects.\footnote{Of course, in order to spell out the truth definition for the language in which the boolean-valued truth definition is stated ($\mathcal{L}_{ZFC}^+$), we would have to move to an even richer metalanguage.} The $\{ 0,1 \}$-valued sets are, on the proposed view, special cases of arbitrary objects: the ``traditional'' universe $V$ is canonically represented in  $V^{(B)}$ \cite[p.~30]{Bell2005}. Thus, I conclude that  the proposed view is logically and conceptually autonomous.

Linnebo and Pettigrew also distinguish \emph{justificational autonomy}: we should be able to justfify the proposed view without appealing to reasons that are specifically tailored to justifying the existence of $\{ 0,1 \}$-valued sets of various kinds \cite[p.~242]{LinneboPettigrew2011}.  I believe that this requirement, too, can be met. Nonetheless, since this article is concerned with articulating and exploring a new foundational proposal rather than arguing that it is superior to other proposals, I leave the question of justificatory autonomy aside here.

On the proposed view, $\mathcal{V}$ is the ``ultimate'' set theoretic universe. In this sense, an \emph{absolutist} interpretation of set theory is proposed. Nevertheless, there is an obvious connection with multiverse views such as that of \cite{Hamkins2012}, \cite{Steel2014}, \cite{Vaananen2014}. We have seen how every state that $\mathcal{V}$ can be in determines a boolean-valued set theoretic universe $V^{(B)}$. Moreover, if we take an anti-chain $\mathcal{A}$ in $B$, then every $a\in \mathcal{A}$ determines a set-theoretic universe. The universes determined by elements of $\mathcal{A}$ will in general not be classical 2-valued universes: they are boolean-valued universes. Moreover, such universes \emph{themselves} typically contain other universes. In this sense, the ultimate set theoretic universe contains many `multiverses'. So the position under consideration can be labeled \emph{kaleidoscopic} absolutism.\footnote{The multiverses in $V^{(B)}$ (determined by anti-chains) can be turned into multiverses of classical, two-valued universes by well-known ultrafilter techniques \cite[chapter 4]{Bell2005}.}

As a foundational mathematical theory, set theory must be sufficiently rich to carry out all of accepted mathematics, albeit sometimes in an exceedingly cumbersome way. Thus, in a naturalistic spirit, I take it as a \emph{conditio sine qua non} that the set-theoretic universe makes $ZFC$ true, and we have seen that $\mathcal{V}$ does this.

As mentioned above, there is a 2-valued universe $V$ that is canonically embedded in every $V^{(B)}$. But the idea is that our mathematical experience suggests that the set theoretic world is not such a 2-valued structure
\cite[p.~418]{Hamkins2012}:
\begin{quote}
[The] abundance of set-theoretic possibilities poses a serious difficulty for the universe
view, for if one holds that there is a single absolute background concept of set, then one
must explain or explain away as imaginary all of the alternative universes that set theorists
seem to have constructed. This seems a difficult task, for we have a robust experience in
those worlds, and they appear fully set theoretic to us. 
\end{quote}
Hamkins takes the independence phenomena to be evidence for his multiverse view; I take them to be evidence for the kaleidoscopic absolutist view.


One might wonder whether it is reasonable to expect $V$ to be a \emph{state} of (some, or even of every) $V^{(B)}$. If it is, for some $V^{(B)}$,  then $B$ will have at least one \emph{atom} $a$, and $V = V^{(B)}_a$. So then $V^{(B)}$, and therefore also $\mathcal{V}$, will contain at least one maximally specific state, i.e., a possible world in the Carnapian sense of the word.
There are, however, reasons for believing that $V$  is \emph{not} a state of a $V^{(B)}$ in the range of $\mathcal{V}$.\label{constraints} If for some $a\in B$, $V = V^{(B)}_a$, then there is at least one completely classical state that the set theoretic universe can be in. Moreover, this state is then also the \emph{only} fully determinate state that the set theoretic universe can be in. Set theoretic experience provides no reason to think that there is any such super-special universe that the set theoretic universe can be. 

The general picture is then as follows. Set theoretic experience---forcing, combined with large cardinal axioms, infinitary combinatorics,$\hdots$---sugg\-ests that there are \emph{many} states that the set theoretic universe can be in. So $\mathcal{V}$ has to be such that it can be in all and only those states. And this imposes restrictions on what $V$ is like and what $B$ can be like. We have seen that for any $ V^{(B)}$ in the value range of $\mathcal{V}$, the boolean algebra $B$ should probably not be the $\{ 0,1 \}$-algebra. 

Beyond this, matters are less clear.
Since we want universes with at least some large cardinals to be possibilities, we probably do not want $V$ to be G\"{o}del's constructive universe $\mathrm{L}$. Perhaps it can be argued that $V$ contains many large cardinals. If indeed $ V \neq\mathrm{L}$, then elementary considerations concerning forcing show that $\mathcal{V}$ cannot be in the state of being $\mathrm{L}$ (rather than $\mathrm{L}$ just being \emph{definable} in some such state). I take this to be in agreement with the fact that the existence of large cardinals is much less controversial than, for instance, the assertion that the Continuum Hypothesis is true (or the assertion that it is false). The fact that we have only very limited knowledge of what the $ V^{(B)}$'s in the range of $\mathcal{V}$ are, does not, however, preclude us from drawing some conclusions about truth value determinateness beyond the theorems of ZFC. For instance, since forcing techniques show that CH is independent even of set theory with large cardinals, we have good reasons to believe that CH does not have a determinate truth value.



\section{Identity}

In  boolean-valued models we can have $[\![ \xi = \eta ]\!] = a$ for some boolean value $0<a<1$. So it seems that we are committing ourselves to identity being to some degree an \emph{indeterminate} relation.

Evans held that indeterminacy of identity is incoherent \cite{Evans1978}. His argument is a simple \emph{reductio} based on Leibniz' principle of the indiscernibility of identicals. Consider any $\xi$ and $\eta$ that are not determinately identical. Then $\xi$ has a property, viz. being identical to $\xi$, that $\eta$ does not have. So $\xi$ and $\eta$ are determinately different from each other.

It has been observed that, strictly speaking, Evans' argument does not go through. In Evans' argument, Leibniz' principle is applied to the predicate $\lambda z [ z \textrm{ is identical with } \xi]$. But then we can only conclude that $\xi$ and $\eta$ are not identical, not that they are \emph{determinately} non-identical.

However, Williamson (\cite[section 8]{Williamson2005}) has shown how an Evans-like argument nonetheless can be carried through with the use of two further plausible principles. First, the following inference rule seems valid:
\begin{quote}
``From a proof of $\phi \rightarrow \psi$, infer that if it is determinately the case that $\phi$, then it is determinately the case that $\psi$.''
\end{quote}
Moreover, if it is determinately the case that $\phi$, then $\psi$. Using these proof principles, Evans' argument can validly be strengthened to conclude that there can be no $\xi$, $\eta$, such that (1) it is determinately the case that they are not determinately identical and (2) it is also determinately the case that they are not determinately different. 

The moral that is often taken from arguments such as these is that there is no ontological vagueness but only semantic vagueness. That is, I surmise, also the attitude that set theorists habitually take, and it is perhaps the main reason why the $V^{(B)}$'s other than $V$ are not taken to be candidates for being the `real' mathematical universe. Indeed, in the forcing poset approach, the vagueness involved is pretty much \emph{officially} regarded as semantic in nature, for its counterparts of the `ontologically vague' sets in the boolean-valued approach are the $\mathbb{P}$\emph{-names} \cite[chapter 7, \textsection 2]{Kunen1980}. 

I believe that the received view that there is only semantic (and epistemic) vagueness, is correct. I will now argue that this view is compatible with the foundational proposal that is explored in the present article.

Clearly there are many pairs $\xi$, $\eta$ of boolean-valued sets that are not numerically identical to each other but that are ``judged'' to be identical by certain boolean-valued models. Here is a simple example:

\begin{example} 
Consider the simple boolean algebra $B_0 = \{ 0,a,b,1   \} $ with $0<a,b<1$ and $a\bot b$. Let $\xi = \{ \emptyset \rightarrow 1 , (\emptyset \rightarrow 1)\rightarrow 1  \}, $ meaning that $dom(\xi) = \{  \emptyset, \emptyset \rightarrow 1 \} $ and $\xi(\emptyset) = \xi(\emptyset \rightarrow 1) = 1.$ Moreover, let 
$\eta = \{ u \rightarrow 1, v\rightarrow 1  \},$ with $u$ and $v$ being the following ``anti-correlated'' sets:
\begin{itemize}
\item $u= \{ \emptyset \rightarrow a, (\emptyset \rightarrow 1) \rightarrow b \}$;
\item $v= \{ \emptyset \rightarrow b, (\emptyset \rightarrow 1) \rightarrow a \}$.
\end{itemize}
Then clearly we have, in the strict sense, $\xi \neq \eta$. Nonetheless, a routine but tedious calculation shows that $V^{(B_0)}\models \xi = \eta$.
\end{example}

In a boolean-valued class model, the identity symbol `=' expresses a congruence relation other than the \emph{real} identity relation: it `measures' the states in which its arguments \emph{coincide}.
On the proposal under consideration, some boolean-valued class models are good interpretations of the language of set theory. Therefore the proposed view is committed to the claim that in set theory, the symbol `=' does not express the real identity relation. As a consequence, it is not threatened by Evans' argument, nor by Williamson's modification of it.

It is at the metaphysical level, i.e., in arbitrary object theory, that we truthfully say that $\xi \neq \eta$; in set theory, we truthfully say that $\xi = \eta$. So in these two contexts we do not use the identity symbol with the same meaning. From  the debate about mathematical structuralism we are familiar with the claim that in many areas of mathematics the identity symbol is commonly not used to express the metaphysical relation of identity---remember the slogan``identity is isomorphism''. Set theory, as a foundational discipline, is often taken to be an exception to this phenomenon. On the view that is explored in the present article, this is not correct: the identity symbol in set theory also expresses a relation that is different from ``real'' identity.



\section{Closing remarks}

On the view that I have sketched, \emph{there is a set theoretic universe}. Moreover, as an arbitrary object, \emph{it is an abstract entity}. It seems natural to say that on the proposed conception the set theoretic universe is \emph{mind-independent}. 
The combination of these three commitments makes the position under consideration a form of \emph{mathematical platonism}. At the same time, this position rejects a strong form of \emph{truth value realism} according to which every set-theoretic statement has exactly one of the traditional truth values (true, false). However, it is now fairly generally recognised that mathematical platonism \emph{per se} is not committed to this extra thesis, even though most traditional forms of mathematical platonism do sign up to it.

Versions of set theoretical platonism without truth value realism have been proposed in the literature. In this article, I have suggested one particular such view that takes the set theoretic universe and the sets in it to be arbitrary objects.
I do not claim that the view that I have proposed is the \emph{only} way in which forcing models can metaphysically be related to arbitrary object theory. I will close by outlining the contours of an alternative way of seeing elements of forcing models as arbitrary objects.

The construction of a forcing model is sometimes seen as analogous to the process of adjoining an object to an algebraic structure.\footnote{See \cite[p.~2]{Chow2008}, for instance.} For definiteness, consider the construction of the ring of polynomials in one variable over $\mathbb{R}$. In terms of arbitrary object theory, this process can roughly be seen as follows. We start with an arbitrary object $X$ with value range $\mathbb{R}$. Then we consider all arbitrary objects that \emph{depend} on $X$ in the sense of being polynomially determined by $X$. The resulting collection of arbitrary objects form a ring.

Similarly, given a poset $\mathbb{P}$ in a model of set theory $M$, a generic filter $G$ can be taken to be an arbitrary subset of $\mathbb{P}$. This is the view of Venturi, who motivates it as follows \cite[p.~2]{Venturiforthc}:
\begin{quote}
Intuitively a set is generic, with respect to a model $M$ and a poset $\mathbb{P}$, if
it meets all requirements to be a subset of $\mathbb{P}$ from the perspective of $M$ and
nothing more. The elements of $\mathbb{P}$, called \emph{conditions}, represent partial pieces
of information that will eventually give the full description of the generic
$G$. Moreover, the dense sets that belong to $M$ represent the properties that
a subset of $\mathbb{P}$ should eventually have, as considered from the perspective of
$M$. For this reason, a generic set does not have a characteristic property
that distinguishes it from all other elements of $M$.
\end{quote}
In terms of this arbitrary subset $G$ of $\mathbb{P}$, a model $M[G]$ can then be seen as a collection of \emph{dependent} arbitrary objects.

This way of connecting forcing models with arbitrary object theory makes use of a distinction between \emph{dependent} and \emph{independent} arbitrary objects. Such a distinction does not figure in the theory of arbitrary objects that I favour.\footnote{Not fundamentally, anyway; but notions of dependence can be \emph{defined} in the framework that was used in the present paper: see \cite[section 9.4]{Horsten2019}.} But it is the cornerstone of Fine's arbitrary object theory. So the alternative account that I have tried to outline in this section is perhaps best developed fully within the framework of Fine's theory. A comparison of this alternative account with the view that was the focal point of this article is left for another occasion.

\end{document}